\documentstyle[12pt, epsf, epsfig, amssymb, amstext, amstex,amsthm, amsfonts]{article}

\input amssymb.sty

\font\tenmsb=msbm10
\font\sevenmsb=msbm7
\font\fivemsb=msbm5
 
\newfam\msbfam 
\textfont\msbfam=\tenmsb
\scriptfont\msbfam=\sevenmsb
\scriptscriptfont\msbfam=\fivemsb
 
\def\Bbb#1{\fam\msbfam\relax#1}

\newtheorem{thm}{Theorem}[section]
\newtheorem{prop}[thm]{Proposition}
\newtheorem{cor}[thm]{Corollary}
\newtheorem{lem}[thm]{Lemma}
\newtheorem{conj}[thm]{Conjecture}
\newtheorem{exa}[thm]{Example}
\newtheorem{defn}[thm]{Definition}

\newtheorem{rem}[thm]{Remark}
\newtheorem{note}[thm]{Notation}
\newtheorem{alg}[thm]{Algorithm}

\newcommand{\ben}{\begin{enumerate}}
\newcommand{\een}{\end{enumerate}}
\newcommand{\ble}{\begin{lem}}
\newcommand{\ele}{\end{lem}}
\newcommand{\bth}{\begin{thm}}
\renewcommand{\eth}{\end{thm}}
\newcommand{\bpr}{\begin{prop}}
\newcommand{\epr}{\end{prop}}
\newcommand{\bco}{\begin{cor}}
\newcommand{\eco}{\end{cor}}
\newcommand{\bcon}{\begin{conj}}
\newcommand{\econ}{\end{conj}}
\newcommand{\bde}{\begin{defn}}
\newcommand{\ede}{\end{defn}}
\newcommand{\bex}{\begin{exa}}
\newcommand{\eex}{\end{exa}}
\newcommand{\brem}{\begin{rem}}
\newcommand{\erem}{\end{rem}}
\newcommand{\bnot}{\begin{note}}
\newcommand{\enot}{\end{note}}
\newcommand{\balg}{\begin{alg}}
\newcommand{\ealg}{\end{alg}}

\newcommand{\bib}{thebibliography}

\newcommand{\C}{{\Bbb C}}
\newcommand{\N}{{\Bbb N}}
\newcommand{\PP}{{\Bbb P}}
\newcommand{\R}{{\Bbb R}}

\begin{document}

\newcommand{\fnref}[1]{~(\ref{#1})}
\newenvironment{emphit}{\begin{itemize} \em}{\end{itemize}}
\begin{center}
\Large{\bf {Identifying Powers of Half-Twists and Computing its Root}}\\ 
\vspace{7mm}
\large{T. Ben-Itzhak, 
S. Kaplan
 and M. Teicher}
\footnote{Partially supported by the Emmy Noether Research Institute for Mathematics, the Minerva 
Foundation of Germany, the Excellency Center "Group Theoretic Methods in the Study of Algebraic 
Varieties"  of the Israel Science Foundation, and by EAGER (European Network in Algebraic 
Geometry).\\}
\end{center}

\begin{abstract}
In this paper we give an algorithm for solving a main case of the conjugacy
problem in the braid groups. We also prove that half-twists satisfy a
special root property which allows us to reduce the solution for the
conjugacy problem in half-twists into the free group. Using this algorithm
one is able to check conjugacy of a given braid to one of E. Artin's generators in any power, and 
compute its root. Moreover, the braid element which conjugates a given
half-twist to one of E. Artin's generators in any power can be restored.
The result is applicable to calculations of braid monodromy of branch curves and verification of 
Hurwitz equivalence of braid monodromy factorizations, which are essential in order to determine 
braid monodromy type of algebraic surfaces and symplectic 4-manifolds.
\end{abstract}

\section*{Introduction}

During past decades braid groups have become important in many fields. Hence, a practical solution 
for its conjugacy problem has become extremely important. Although the groups conjugacy problem was 
first solved by Garside \cite{GAR} (1969) and was addressed many times in the past (i.e., 
\cite{EFF}, \cite{NEW}), still a practical polynomial algorithm for its solution is unknown.

This has lead to the research for the solution of partial problems such as identification of 
special conjugacy classes. In \cite{KAPLAN}, a random algorithm for the identifying half-twists in 
any power was given, and the aim of this paper is to give a deterministic algorithm for the 
problem, which although exponential is of interest because of the simplicity of the proofs 
involved, and the combination of techniques used in order to solve the problem.

The algorithm presented here enables to test factors of braid monodromy type of surfaces, and is a 
partial solution for the conjugacy problem in the braid groups. Moreover, it enables to solve 
specific cases of quasi-positivity.

We begin in Section $1$, by giving some basic definitions related with braid groups. Then, in 
Section $2$, we prove a simple manner of the root property conjectured and tested by a computer in 
\cite{KAPLAN} for half-twists. This means that if we know that two identical powers of half-twists 
are conjugated to one another by a positive braid word $w$, then the two half-twists are conjugated 
by the same word $w$. In section $3$ we give some properties of powers of half-twists which will be 
used in Section $4$ to present the algorithm for identifying the powers of half-twists in the braid 
group and to compute an elemenent which conjugates them in any power to one given Artin generator.


\section{Braid group preliminaries}

\subsection{Artin's braid group}

In this section we will give the definition and some properties of the braid group.

\bde \label{Presentation}
\underline{Artin's braid group $B_n$} is the group generated by $\{\sigma _1,...,\sigma _{n-1}\}$
submitted to the relations 
\ben
\item $\sigma _i \sigma _j=\sigma _j \sigma _i$ where $|i-j| \geq 2$ 
\item $\sigma _i \sigma _{i+1}\sigma _i=\sigma _{i+1}\sigma _i \sigma _{i+1}$  for all 
$i=1,...,n-2$
\een
\ede

This algebraic definition can be looked from a geometric point of view, by associating to every 
generator
of the braid group $\sigma _i$ a tie between $n$ strings going monotonically from top to bottom, 
such
that we switch by a positive rotation between the two adjacent pair of strings $i$ and $i+1$.
The operation for the geometric group is the concatenation of two geometric 
sets of strings. 

\begin{figure}[h]
\begin{center}
\epsfxsize=2cm
\epsfbox{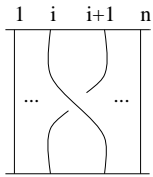}
\caption{The geometrical braid associated with $\sigma _i$.}
\end{center}
\end{figure}

\bde
There is a natural homomorphism, denoted by $\pi$ from the braid group to the symmetric group, 
sending each braid to the permutation induced by the braid strings. $\pi$ is defined by $\pi 
(\sigma _i) = (i,i+1)$
\ede

\noindent
We need some definitions:

\bde
Let $w \in B_n$ be a braid. Then it is clear that $w=\sigma _{i_1}^{e_1} \cdot ...\cdot \sigma
_{i_l}^{e_l}$ for some sequence of generators, where $i_1,...,i_l \in \{1,...,n-1\}$ and 
$e_1,...,e_l \in \{1,-1\}$. We will call such presentation of $w$ a \underline{braid
word}, and $\sigma _{i_k}^{e_k}$ will be called the \underline{$k^{th}$ letter of the word $w$}. 
$l$ is the \underline{length} of the braid word, and we define $exp(w) = \Sigma _{i=1} ^l e_i$, 
which is the sum of exponents of the letters of $w$, and known to be invariant in the conjugacy 
class of $w$.
\ede

We will distinguish between two relations on the braid words.

\bde
Let $w_1$ and $w_2$ be two braid words. We will say that $w_1=w_2$ if they represent the
same element of the braid group.
\ede

\bde
Let $w_1$ and $w_2$ be two braid words. We will say that $w_1 \equiv w_2$ if $w_1$ and
$w_2$ are identical letter by letter.\
\ede

\bde
A positive braid is an element of $B_n$ which can be written as a word in positive powers of the
generators $\{\sigma _i\}$, without the use of the inverse elements $\sigma _i^{-1}$. We denote
this subsemigroup $B_n^+$
\ede


\bde
Let $\beta \in B_n$, consider two different strings, the i-th and the j-th $i \neq j$. In a general 
diagram of $\beta $ these two strings ``intesect'' at several points. Let $p(i,j)$ (and $n(i,j)$) 
denote the number of times the i-th string crosses under the j-th string starting from the left 
(from the right). Then the number $p(i,j) - n(i,j)$ is called the crossing index from the i-th 
string to the j-th string denoted by $cr_{\beta} (i,j)$.
\ede

\noindent The following Lemma is well known:

\ble \label{CRLemma1}
\ele
\noindent
1. The crossing index of two strings is an invariant of the braid.\\
2. If $\pi (\beta ) = Id$ then $cr_{\beta } (i,j) = cr_{\beta } (j,i) \quad \forall i,j$.\\
3. $\Sigma _{i \neq j} cr _{\beta } (i,j) = exp(\beta )$.


\subsection{The half-twists}

We are going to describe a specific conjugacy class in the braid groups which we will call its 
elements Half-twists.
 
\bde \label{ConjDef}
Let $H$ be the conjugacy class of $\sigma _1$, (i.e. $H=\{q^{-1} \sigma _1 q : q \in
B_n$\}). We call $H$ \underline{the set of half-twists in $B_n$}, and we call an element $\beta 
\in H$ a \underline{half-twist}.
\ede

\ble \label{Transposition}
If $b \in H$ is a half-twist then there exists $1 \leq i,j \leq n$, such that $\pi (b)=(i,j)$.
\ele
\noindent {\bf Proof: } Trivial, since $\pi(P \sigma _i P^{-1}) = \pi(P)^{-1} (i,i+1) \pi(P)$ which 
is also a transposition.\\

\noindent It is proven in \cite{BGTI} that the half-twists occupy a full conjugacy class in the 
braid group. Moreover, it is proven in \cite{BGTI} that all the generators of the braid group are 
conjugated to each other.

\subsection{Half-twists and paths}
Another way of looking at the braid group on $n$ generators is by the group of isotopy classes of 
orientation preserving automorphisms of an $n$ punctured disk that fix the set of punctures and is 
the identity on the boundary of the disk. It is proven in \cite{BGTI} that the all half-twists can 
be regarded as an automorphism that is described by a path connecting two punctures which does not 
self intersect nor intersect the punctures aside for its beginning and ending points. The 
automorphism induced by this path actually exchange the two punctures clockwise in a small 
neighborhood of the it, which does not contain any of the other punctures. 

The geometric braid can be determined by taking the trace of the punctures under the isotopy as can 
be seen in the figure below:

\begin{figure}[h]
\begin{center}
\epsfxsize=3.5cm
\epsfbox{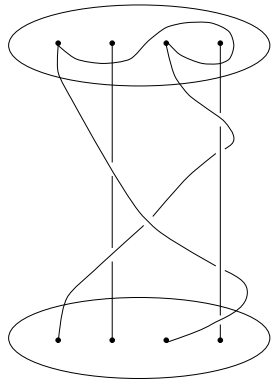}
\label{Fig2}
\caption{The geometric braid and its defining path.}
\end{center}
\end{figure}

\brem \label{SurfaceRemark}
This actually means that for half-twist braids, we have two strings which interact along a small 
neighborhood of a surface in $\R ^3$, while the other strings go simply straight from top to 
bottom. The surface is what is found underneath the path.
\erem 

\section{The root property of half-twists}

In \cite{KAPLAN}, we conjectured after testing with a computer, that a braid word $w$ is conjugated 
to the generator $\sigma _i$ if and only if $w^2$ is conjugated to $\sigma _i^2$.

The aim of this section is to give a proof for theorem \ref{THE_THEOREM} which is a simple version 
of the above conjecture, and says that if $W \sigma _i^m =\sigma _j^m W$ for some $W \in B_n^+$, $m 
\in \N$. Then $W \sigma _i =\sigma _j W$.
However, before we can go to the proof we need some properties of positive braids.

\subsection{Properties of positive braids}

In this section we will recall two properties of positive braids. The first is concerning with the 
equivalence 
of two positive braid words, and the second is concerning the ability to write every conjugated 
element in the braid group using 
a positive braid word as its conjugator.

\bpr 
If $\alpha=\beta \in B_n^+$, then the length $|\alpha |$ of $\alpha$ is equal to the length $|\beta 
|$ of $\beta$.
\epr

\noindent
{\it Proof:} This follows immediately from the relations of the braid group since they do not 
change the number of generators.
\hfill $\qed$

\bde
Let $\alpha,\beta \in B_n^+$ be two positive braid words. We call $\alpha$ and $\beta$ positively 
equivalent if there is a sequence of positive words $\alpha \equiv w_0 = w_1 = ... = w_k \equiv 
\beta$, such that $w_{i+1}$ is obtained from $w_i$ by a single activation of a relations of $B_n$.
\ede

\bpr
Let $\alpha,\beta \in B_n^+$ be two braid words such that $\alpha=\beta$ in $B_n$. Then $\alpha$ 
and $\beta$ are positively equivalent. 
\epr

\begin{proof}
See \cite{GAR}
\end{proof}

\bpr \label{POS_ALG}
Let $\alpha=\beta \in B_n^+$, such that $\alpha \equiv \sigma _i w$ for some generator $\sigma _i$ 
and a positive braid word $w$. 
Then, there is an algorithm that writes $\beta \equiv \sigma _i w'$ for some positive braid word 
$w'$ such that $|w|=|w'|$. 
\epr 

\noindent
{\it Proof:}
The proof can be found in \cite{EFF}.
\hfill $\qed$

If we examine this algorithm in detail we result with the following corollary: 

\bco \label{co1}
Let $\alpha \equiv \sigma _i \sigma _j w \in B_n^+$ for some generators $\sigma _i,\sigma _j$ $i 
\neq j$ and a positive braid 
word $w$, and suppose that $\alpha=\sigma _j w'$ for some $w' \in B_n^+$. Then, one of the 
following must happen:
\ben
\item
If $|i-j|>1$ we can write $\alpha=\sigma _j \sigma _i w$.

\item
If $|i-j|=1$ we can write $w=\sigma _i w''$
\een
\eco

\noindent
{\it Proof:}
Although the proof is obvious from the proof given in \cite{EFF} for the correctness of the 
algorithm in Proposition \ref{POS_ALG} that extract a letter to the 
left of a braid word, we will give the basic idea of it. In the first case simply $\sigma _i$ and 
$\sigma _j$ commute, so write 
$w=w'$ and we have finished. In the second case the only way to transform to the left the letter 
$\sigma _j$ is to use the second 
braid relation $\sigma _i \sigma _{i+1} \sigma _i=\sigma _{i+1} \sigma _i \sigma _{i+1}$. Now since 
it is given that 
$\alpha=\sigma _j w'$ for some $w' \in B_n^+$, we must be able to activate the braid relations in 
order make this switch, 
but this means that we can write $w=\sigma _i w''$ for some $w'' \in B_n^+$.
\hfill 
$\qed$

\bco \label{co-2.4}
Let $\alpha \equiv \sigma _i^k \sigma _j w \in B_n^+$ for some generators $\sigma _i,\sigma _j$ $i 
\neq j$ and a positive braid 
word $w$, and suppose that $\alpha=\sigma _j w'$ for some $w' \in B_n^+$. Then, one of the 
following must happen:
\ben
\item
If $|i-j|>1$ we can write $\alpha=\sigma _j \sigma _i^k w$.

\item
If $|i-j|=1$ we can write $w=\sigma _i w''$
\een
\eco

\noindent
{\it Proof:}
The proof is clear from the proof of Corollary \ref{co1} and the algorithms in \cite{EFF}.

\hfill
$\qed$

\bco \label{co2}
If we use the second case in the corollary above then $|w''|<|w|$. 
\eco 

The following theorem will be of great importance in what follows.

\bth \label{th1}
Let $\alpha ,\beta \in B_n$ (not necessarily positive) two conjugated braids. Then, there exist a 
positive braid $w \in B_n^+$ 
such that $w^{-1} \alpha w = \beta$.
\eth

\begin{proof}
See \cite{POSBR}
\end{proof}

\subsection{Proof of The Theorem}

Before we start with the theorem and its proof we want 
to point out that the unique root property, that we are going to prove for the conjugacy class of 
the half-twists, does not hold for the entire braid group elements. For example:

\bex
Consider the two braids in $B_4$:

$$B=\sigma _1 \sigma _2 \sigma _3 \sigma _1 \sigma _2 \sigma _3 $$
$$\Delta=\sigma _1 \sigma _2 \sigma _3 \sigma _1 \sigma _2 \sigma _1 $$
\eex

\bpr
The two braids described above are different.
\epr

\noindent
{\it Proof:} If $B=\Delta$ then, 

$$B=\sigma _1 \sigma _2 \sigma _3 \sigma _1 \sigma _2 \sigma _3 
=\sigma _1 \sigma _2 \sigma _3 \sigma _1 \sigma _2 \sigma _1=\Delta$$
$$\sigma _3=\sigma _1$$
a contrediction.

\hfill $\qed$

But, it is a common knowledge that the square of both $B$ and $\Delta$ equals $\Delta ^2$ which is 
the generator of the cyclic 
center of the braid group.

\vspace{12pt}
\noindent
The rest of this section is dedicated to the proof of the following theorem:

\bth \label{RootTheorem}
Let $H_1$,$H_2$ be two half-twists such that $H_1^m=H_2^m$. Then, $H_1=H_2$.
\eth 

This means that inside the conjugacy class of the half-twists the braid group has the unique root 
property. 

First, we want to reduce the problem into an equivalent one. Since we know that $H_1$ and $H_2$ are 
two half-twists we know that 
there are two positive braid words $w_1,w_2 \in B_n^+$ such that $H_1=w_1^{-1} \sigma _i w_1$ and 
$H_2=w_2 \sigma _j w_2^{-1}$.
Therefore, since $H_1^m=H_2^m$ we have that $w_1^{-1} \sigma _i^m w_1=w_2 \sigma _j^m w_2^{-1}$. 
But this means that 
$w_2^{-1}w_1^{-1} \sigma _i^m w_1w_2=\sigma _j^m$.
If we can prove that in this case $w_2^{-1}w_1^{-1} \sigma _i w_1w_2=\sigma _j$. we result 
immidietly with what we want.
Therefore, we can reduce the theorem into the following theorem:

\bth \label{THE_THEOREM}
Let $W \sigma _i^m =\sigma _j^m W$ for some $W \in B_n^+$, $m \in \N$. Then $W \sigma _i =\sigma _j 
W$.
\eth

\noindent
{\it Proof:}
We will proof this theorem by induction on the length $|W|$ of $W$.

\underline{The base}:

\ben
\item
If $|W|=0$ it means that $\sigma _i^m=\sigma _j^m$. Threrfore $i=j$ and we have that $\sigma 
_i=\sigma _j$ as we need.
\item
If $|W|=1$ we have $\sigma _k \sigma _i^m=\sigma _j^m \sigma _k$. 

If $|k-i|>1$ or $|k-j|>1$ or $k=i$ or $k=j$ then, at least 
one of $\sigma _i$ or $\sigma _j$ commutes with $\sigma _k$, and therefore we have (without loss of 
generality) 
$\sigma _k \sigma _i^m=\sigma _k \sigma _j^m$, which means that $\sigma _i^m=\sigma _j^m$.

If $|i-k|=1$, we have $2$ different cases:
\ben
\item
If $k=i+1$ and $k=j+1$, or $k=i-1$ and $k=j-1$ we have $i=j$, but then when checking the linking 
numbers associated with, $\sigma _k \sigma _i^m \sigma _k^{-1}$ and $\sigma _j^m=\sigma _i^m$ we 
find a contradiction, which means that this case can not happen. 
\item
If $k=i+1$ and $k=j-1$, or $k=i-1$ and $k=j+1$ we have again that the linking numbers associated 
with the braids $\sigma _k \sigma _i^m \sigma _k^{-1}$ and $\sigma _j^m$ give us a contradiction, 
which means that this case can not happen.
\een 
\een

Now, suppose that for all $W$ such that $|W|<n$ the following is true $W \sigma _i^m=\sigma _j^m W 
\Rightarrow W \sigma _i=\sigma _j W$.
We need to prove that if $|W|=n$ then, $W \sigma _i^m=\sigma _j^m W \Rightarrow W \sigma _i=\sigma 
_j W$.

Now lets look at $W$. We divide the proof into some cases:

\ben
\item
We can write $W=\sigma _k W'$ where $W' \in B_n^+$ and $|k-j|>1$ or $k=j$. 

Then, 
$$\sigma _k W' \sigma _i^m=W \sigma _i^m=\sigma _j^m W=\sigma _j^m \sigma _k W'=\sigma _k \sigma 
_j^m W'$$
$$\Downarrow$$
$$W' \sigma _i^m=\sigma _j^m W'$$
and we finished by induction since $|W'|<|W|$.

\item
We can write $W=W' \sigma _k$ where $W' \in B_n^+$ and $|k-i|>1$ or $k=i$. Then, 
$$W' \sigma _i^m \sigma _j=W' \sigma _k \sigma _i^m=W \sigma _i^m=\sigma _j^m W=\sigma _j^m W' 
\sigma _k$$
$$\Downarrow$$
$$W' \sigma _i^m=\sigma _j^m W'$$
and we finished by induction since $|W'|<|W|$.

\item
If none of the above cases is true, it means that we can't write $W$ beginning or ending with a 
letter that commutes with 
$\sigma _i$ or with $\sigma _j$. So we know that in every way of writing $W=\sigma _kW'\sigma _l$ 
where $W' \in B_n^+$ and 
$|k-j|=1$ and $|l-i|=1$. Now we have to split again into cases:
  
\ben
\item
Suppose that $W=\sigma _{j+1}W'\sigma _{i+1}$
Then by what is given we know that 
\begin{equation} \label{aaa}
W \sigma _i^m=\sigma _{j+1}W' \sigma _{i+1} \sigma _i^m=\underbrace{\sigma _j^m \sigma _{j+1} 
\underbrace{W' \sigma _{i+1}}_{w}}_{\alpha}=\sigma _j^m W
\end{equation} 
Therefore, by \ref{co-2.4} when we take $\alpha=\sigma _j^m \sigma _{j+1} W' \sigma _{i+1}$, and 
$w=W'\sigma _{i+1}$ we have that $w=W'\sigma _{i+1}=\sigma _j W''$.

Substituting this in (\ref{aaa}) we get 
$$\sigma _{j+1} \sigma _j W'' \sigma _i^m=\sigma _j^m \sigma _{j+1} \sigma _j W''=$$
$$=\sigma _j^{m-1} \sigma _{j+1} \sigma _j \sigma _{j+1} W''=\cdots =\sigma _{j+1} \sigma _j \sigma 
_{j+1}^m W''$$
$$\Downarrow$$
$$W'' \sigma _i^m=\sigma _{j+1}^m W''$$

Because of \ref{co2} we know that $|W''|<|W|$ so, by induction we have:

$$W'' \sigma _i=\sigma _{j+1} W''$$

Now, by multiplying this equation by $\sigma _{j+1} \sigma _j$ on the left we get:

$$\sigma _{j+1} \sigma _j W'' \sigma _i=\sigma _{j+1} \sigma _j \sigma _{j+1} W''$$
$$\Downarrow$$
$$\sigma _{j+1} \underbrace{\sigma _j W''}_{W' \sigma _{i+1}} \sigma _i=\sigma _j \sigma _{j+1} 
\underbrace{\sigma _j W''}_{W' \sigma _{i+1}}$$
$$\Downarrow$$
$$W \sigma _i=\sigma _{j+1}W' \sigma _{i+1} \sigma _i=\sigma _j \sigma _{j+1} W' \sigma 
_{i+1}=\sigma _j W$$

\item
A similar argument implies that When $W=\sigma _{j-1}W'\sigma _{i-1}$ the theorem holds.

\item
Suppose that $W=\sigma _{j-1}W'\sigma _{i+1}$
Then by what is given we know that 
\begin{equation} \label{bbb}
W \sigma _i^m=\sigma _{j-1}W' \sigma _{i+1} \sigma _i^m=\underbrace{\sigma _j^m \sigma _{j-1} 
\underbrace{W' \sigma _{i+1}}_{w}}_{\alpha}=\sigma _j^m W
\end{equation}

Therefore, by \ref{co-2.4} when we take $\alpha=\sigma _j^m \sigma _{j-1} W' \sigma _{i+1}$ and 
$w=W' \sigma _{i+1}$ we must be able to write $w=W' \sigma _{i+1}=\sigma _j W''$.

Substituting this in (\ref{bbb}) we get
$$\sigma _{j-1} \sigma _j W'' \sigma _i^m=\sigma _j^m \sigma _{j-1} \sigma _j W''$$
However, $\sigma _j^m \sigma _{j-1} \sigma _j W''=\sigma _j^{m-1} \sigma _{j-1} \sigma _j \sigma 
_{j-1} W''=\cdots =\sigma _{j-1} \sigma _j \sigma _{j-1}^m W''$, hence we have:

$$\sigma _{j-1} \sigma _j W'' \sigma _i^m=\sigma _{j-1} \sigma _j \sigma _{j-1}^m W''$$
$$\Downarrow$$
$$W'' \sigma _i^m=\sigma _{j-1}^m W''$$

Because of \ref{co2} we know that $|W''|<|W|$ so, by induction we have:

$$W'' \sigma _i=\sigma _{j-1} W''$$

Now, multiplying this equation by $\sigma _{j-1} \sigma _j$ on the left we get:

$$\sigma _{j-1} \sigma _j W'' \sigma _i=\sigma _{j-1} \sigma _j \sigma _{j-1} W''$$
$$\Downarrow$$
$$\sigma _{j-1} \underbrace{\sigma _j W''}_{W' \sigma _{i+1}} \sigma _i=\sigma _j \sigma _{j-1} 
\underbrace{\sigma _j W''}_{W' \sigma _{i+1}}$$
$$\Downarrow$$
$$W \sigma _i=\underbrace{\sigma _{j-1} W' \sigma _{i+1}}_{W} \sigma _i=\sigma _j 
\underbrace{\sigma _{j-1} W' \sigma _{i+1}}_{W}=\sigma _j W$$

\item
A similar argument shows that the theorem holds when $W=\sigma _{j+1}W'\sigma _{i-1}$
\een 
\een

\hfill 
$\qed$

\section{Properties of half-twists powers}

\noindent In the following section we introduce some properties of half-twists and power of 
half-twists. In particular, we prove that the problem of determening if a braid is an event power 
of a half-twist is equivalent to the conjugacy problem in the free group.


\bde
A braid which is equivalent to a braid that can be drawn in a way, where strings $2,...,n$ are 
parallel and straight, is called a combed braid, the process of transforming a given combed braid 
into it's equivalent representation with strings $2,...,n$ straight is called {\it combing}.
\ede


\ble \label{CombedLemma}
\ele
\noindent
1. The set of all combed braids is a subgroup of $B_n$ denoted by $A_n$.\\
2. $A_n$ is the free group generated by $\{ a_1,...,a_{n-1}\}$, where,
$$a_i = \sigma _1... \sigma _{i-1} \sigma _i ^2 \sigma _{i-1} ^{-1}... \sigma _1 ^{-1}$$
{\bf Proof: } \cite{Artin}.\\\\

\noindent The algorithm for getting the presentation of a combed braid in terms of the 
$\{a_1,...,a_{n-1} \}$ generators is well known. We recall here its steps, for clarification. 

\balg \label{Alg1}
\ealg

\begin{figure}[h]
\begin{center}
\epsfxsize=10cm
\epsfysize=6cm
\epsfbox{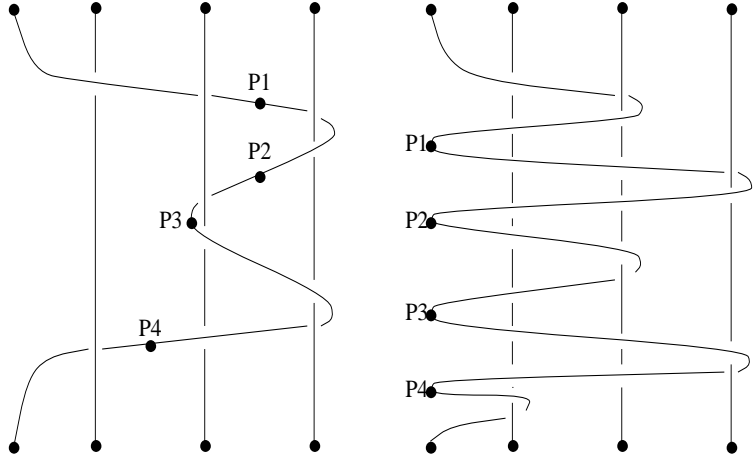}
\label{FigAn}
\caption{Changing to $A_n$ presentation by pulling back all non-matching points.}
\end{center}
\end{figure}

\noindent {\bf Step 1:} 'Comb' the $2,...,n$ strings of the braid so they are all parallel, and 
only the 1st string is moving in between them.\\\\
\noindent {\bf Step 2:} Moving along the first string, every time, the string does not match to the 
$A_n$ presentation, 
'pull' the string at this point to the left side of the braid. An example of Step 2, is shown in 
figure \ref{FigAn}. In the example we show how does the combed braid $\sigma _1 \sigma _2 ^{-1} 
\sigma _3 ^{-2} \sigma _2 ^2 \sigma _3 ^2 \sigma _2 ^{-1} \sigma _1$ can be written in $A_n$ 
presentation as $a_2 ^{-1} a_3 ^{-1} a_2 a_3 a_1 $. \\

It is known that the complexity of Algorithm \ref{Alg1} is exponential. 

\ble \label{CRLemma2}
Let $\beta $ be a half-twist such that $\pi (\beta ) = (i,j)$ and $k > 0$ then $$cr_{\beta 
^{2k}}(s,t) =
{	\begin{cases} 
		{ k \text{ if } \{ s ,t \} = \{ i,j \} } \\
		{ 0 \text{ other } } 
	\end{cases} }
$$
\ele
\noindent {\bf Proof:} Since $\pi (\beta ^2) = Id $, and by Lemma \ref{CRLemma1}, it is clear that 
$$cr_{\beta ^{2k}} (i,j) =  k \cdot cr_{\beta ^2} (i,j)$$
\noindent Therefore, it is enough to prove the lemma for $k=1$.\\\\

\noindent Observing $\beta ^2$ as a trace of $n$ punctures as shown in Figure $2$, we notice the 
following 3 cases:\\\\
\noindent {\bf Case 1:} The s-th and the t-th strings where $\{ s,t \} \bigcap \{ i , j \} = \phi $ 
are parallel, and therefore, $cr_{\beta ^2} (s,t) = 0$.\\
\noindent {\bf Case 2:} Assume that $s=i$ and $t \neq j$ the i-th string crosses the t-th string 
and then moves back in the same path and we get that $p(s,i) = n(i,s)$, therefore $cr_{\beta ^2} 
(i,s)=0$.\\
\noindent {\bf Case 3:} If $\{ s,t \} \bigcap \{ i , j \} =\{ i , j \} $, since, $exp(\beta ^2)$ is 
$2$ and from Lemma \ref{CRLemma1}, 
\begin{equation} \label{eq1}
\Sigma _{s \neq t} cr_{\beta ^2}(s,t) =2
\end{equation}
\noindent which is, from the above, equal to $ cr_{\beta ^2} (i,j) +  cr_{\beta ^2} (j,i)$.\\
Since $ \pi (\beta ^2) = Id$, and from (\ref{eq1}) and Lemma \ref{CRLemma1} (2), $ cr_{\beta ^2} 
(i,j) =  cr_{\beta ^2} (j,i) = 1$.
\hfill $\qed$

\ble \label{FreeConj}
Let $b$ be a half-twist such that $\pi (b) = (1,n)$, then, $b^2$ is combed braid cojugated in $A_n$ 
to $a_{n-1}$
\ele 
\noindent {\bf Proof: } Let $b$ be a half-twist with $\pi (b) = (1,n)$. By definition 
\ref{ConjDef}, $b$ can be written as 
$P^{-1} \sigma _1 P$, for some $P \in B_n$.\\\\
By Remark \ref{SurfaceRemark}, the $1$-st and the $n$-th strings can be embedded in a small 
neighborhood of a surface $F$, which does not intersect any other string.\\
$b^2 = P^{-1} \sigma _1 ^2 P$ and therefore, by changing $\sigma _1$ into $\sigma _1 ^2$ we see 
that $b^2$ can be embedded 
in the same neighborhood of the surface $F$ (See Figure $2$).\\
In $b^2$, strings $1$ and $n$ move along the surface, go around each other and move back along the 
surface to 
their original point. Now, since $F$ does not intersect any other string, we may straight the $n$th 
string along the surface $F$, and make the $n$th string parallel to the 
other $n-2$ parallel strings. Hence, we get that $b^2$ is equivalent to a combed braid (See Figure 
$4$).\\
By the algorithm \ref{Alg1}, getting the $A_n$ generators of the combed braid $b^2$, the braid 
still remain symmetric 
with respect to $\sigma _{n-1} ^2$. Therefore, $b^2$ is conjugated to $a_{n-1}$.

\begin{figure}[h]
\begin{center}
\epsfxsize=12cm
\epsfysize=4cm
\epsfbox{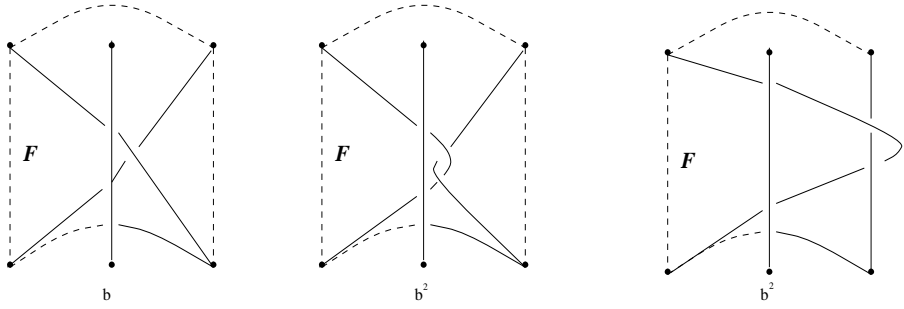}
\label{Fig3}
\caption{$b^2$ is a combed braid}
\end{center}
\end{figure}

\bth \label{CombbedTheorem}
Let $b \in B_n$ with $\pi (b) = (1,n)$, then, $b$ is a half-twist in power of $2k$ if and only if 
$b$ is a combed braid conjugated to $a _{n-1} ^k$ in $A_n$.
\eth
\noindent {\bf Proof: } Let $b = h^{2k}$, where $h$ is a half-twist. By Lemma \ref{FreeConj}, 
$h^{2} = P^{-1} a_{n-1} P$ for some $P \in A_n$ and therefore, 
$$b = h^{2k} = P^{-1} a_{n-1} ^{k} P.$$
\noindent
The second direction of the theorem is trivial since,
$$b = P^{-1} (\sigma _1... \sigma _{n-2} \sigma _{n-1}^2 \sigma _{n-2} ^{-1}... \sigma _1 ^{-1})^k 
P = 
(P^{-1} \sigma _1... \sigma _{n-2} \sigma _{n-1} \sigma _{n-2} ^{-1}... \sigma _1 ^{-1} P)^{2k}$$
which is a half-twist in the power of $2k$.
\hfill $\qed$

\section{The algorithm}

We are now ready to present the algorithm. Given a braid $b$, the algorithm will determine if $b$ 
is a power of a 
half-twist. If the answer is positive, and $b$ is half-twist in the power of $k$, $k \in \N$, the 
algorithm computes a 
$P \in B_n$ such that, $b = P^{-1} \sigma _1^k P$.

\balg \label{THE_ALGORITHM}
\ealg

\noindent
{\bf Step 1:} Compute $exp(b)$. It is clear that if $b$ is a half-twist in power $k$ than $k = 
exp(b)$.\\

\noindent
{\bf Step 2:}\\\\
\noindent
{\bf  Case 1:} If $exp(b)$ is odd, compute $\pi (b)$. If $\pi (b)$ is not a transposition, then, by 
Lemma \ref{Transposition} 
$b$ is not an odd power of a half-twist. Therefore return {\it false}.\\\\
{\bf  Case 2:} If $exp(b)$ is even, compute $\pi (b)$. If $\pi (b) \neq Id$, then, By Lemma 
\ref{Transposition} $b$ is not 
an even power of a half-twist. Therefore, return {\it false}. If $\pi (b) = Id$ compute $cr_b(i,j)$ 
$\forall i,j$ and find $i$ and $j$ such that 
$cr_b(i,j) = cr_b(j,i) = exp(b)/2$ all other $cr_b(k,l)$ must be zero. 
By Lemma \ref{CRLemma2} if such $i$ and $j$ do not exist, or $cr_b (k,l) \neq 0$ where 
$\{i,j\}\neq\{k,l\}$ then $b$ is not a power of a half-twist, and we can return {\it false}.\\

At this point we know, that if $b$ is a power of a half-twist than $b$ is equivalent to a braid 
where all strings other than $i$ and $j$ are parallel. We call the strings $i$ and $j$ the 
switching strings.\\\\

\noindent
{\bf Step 3:} To apply Theorem \ref{CombbedTheorem} we need that the two switching strings found in 
step $2$, $i$ and $j$ will be $1,n$. Assuming that $i<j$, by conjugating 
$\sigma _{i-1}... \sigma _2 \sigma _1$ to $b$ and then conjugating $\sigma _{j+1} ... \sigma 
_{n-1}$ to the result, 
we get a braid $b'$ with this property. Since the half-twists occupy a full conjugacy class $b'$ is 
a power of 
a half-twist if and only if $b$ is a power of half-twist. If we manage to find $P_1 \in B_n$ such 
that, $b' = P_1^{-1}\sigma _1 P_1$
computing $P$ is trivial. At this point of the proof we may assume that the two switching strings 
are $1,n$.\\\\

\noindent
{\bf Step 4:} If $exp(b)$ is odd, Compute $b^2$. If $exp(b)$ is even, $b$ remain the same. At this 
point we know that,
 if $b$ is a power of a half-twist, the result of this step is an even power of a half-twist, and 
therefore a combed braid.\\\\

\noindent
{\bf Step 5:} Operate Algorithm \ref{Alg1} on the result of the previouse step. If the algorithm 
fails (since the braid is 
not a combed braid), By Theorem \ref{CombbedTheorem}, $b$ is not a power of a half-twist, and we 
return {\it false}. On the contrary, if the algorithm succeeded, the 
result of this step is the presentation of $b$ (or $b^2$, depending on the previous step) in 
$A_n$.\\\\

\noindent
{\bf Step 6:} Check if the result is conjugated to $a _{n-1}^{\frac{k}{2}}$ ($a _{n-1}^{k}$ in the 
case of $b^2$) in $A_n$. The conjugacy problem in $A_n$ is easy to compute since 
$A_n$ is a free group (Lemma \ref{CombedLemma}). By Theorem \ref{CombbedTheorem}, the answer is 
positive if and 
only if $b$ (or $b^2$) is a power of a half-twist. So, if the answer is {\it false}, the algorithm 
returns {\it false}.\\\\

\noindent
{\bf Step 7:} In case that the result of step $6$ is true, we compute the k-th (2k-th) root of the 
result: change back to $B_n$ presentation, and since the braid is 
conjugated to $a _{n-1} ^{\frac{k}{2}}$ ($a _{n-1}^{k}$ in the case of $b^2$) we get the form:

$$b = Q^{-1}(\sigma _1... \sigma _{n-2} \sigma _{n-1} ^2 \sigma _{n-2} ^{-1} \sigma _1 
^{-1})^{\frac{k}{2}}Q$$

and the k-th root is:

$$r = Q^{-1}(\sigma _1... \sigma _{n-2} \sigma _{n-1} \sigma _{n-2} ^{-1} \sigma _1 ^{-1})Q$$
which is obviously a half-twist. 

In case of $b^2$ we apply the same procedure.\\
\noindent If $exp(b)$ is even, the algorithm is finished.\\
\noindent If $exp(b)$ is odd, at this point we know that $b^2$ is a half-twist of the power $2k$, 
and we computed its 
2k-th root, $r$. Since the 2k-th root is unique (Theoerm \ref{RootTheorem}), if $b$ is a half-twist 
in the power of k, 
$r$ must be it's k-th root. Therefore, 
in order to check if $b$ is a power of a half-twist it is enough to check whether, $b$ is 
equivalent to $r^k$ in $B_n$. 
There are many solutions for the word problem in the braid group, one can be found for example in, 
\cite{NEW}.\\\\
\noindent This completes the algorithm.\\

\noindent
The complexity of Algorithm \ref{THE_ALGORITHM} is exponential, since in Step $5$ we comb the braid 
which might take exponential time as described after Algorithm \ref{Alg1}. All other steps are 
bounded by this exponential function and hence it gives the upper bound for the number of 
operations which the algorithm performs.



\begin{\bib}{10}
\bibitem{Artin} Artin, E., {\it Theory of braids}, Ann. Math. {\bf 48} (1947), 101-126.
\bibitem{NEW} Birman, J.S., Ko, K.H. and Lee, S.J., {\it A new approach to the word and conjugacy 
problems in the braid groups}, Adv. Math. {\bf 139} (1998), 322-353.
\bibitem{POSBR} Elrifai, E.A. and Morton, H.R., {\it Algorithms for positive braids}, Quart. J. 
Math. Oxford Ser. (2) {\bf 45} (1994), 479-497.
\bibitem{GAR} Garside, F.A., {\it The braid group and other groups}, Quart. J. Math. Oxford Ser. 
(2) {\bf 78} (1969), 235-254.
\bibitem{EFF} Jacquemard, A., {\it About the effective classification of conjugacy classes of 
braids}, J. Pure. Appl. Alg. {\bf 63} (1990), 161-169.
\bibitem{KAPLAN} Kaplan, S. and Teicher, M., {\it Identifying Half-Twists Using Randomized 
Algorithm Methods}, preprint.
\bibitem{BGTI} Moishezon, B. and Teicher, M., {\it Braid group techniques in complex geometry I, 
Line arrangements in $\C \PP^2$}, Contemporary Math. {\bf 78} (1988), 425-555.
\end{\bib}

\end{document}